\documentclass[regno,12pt]{amsart}
\usepackage{amsmath,amsthm,amssymb,a4,graphics,color}
\usepackage[active]{srcltx}
\usepackage[usenames,dvipsnames]{pstricks}
\usepackage{epsfig}
\usepackage{pst-grad}
\usepackage{pst-plot} 
\usepackage{tikz}
\usepackage{hyperref}
\usepackage{dsfont}
\usepackage{float}
\usepackage{todonotes}

\newfam\Bbbfam 
\font\tenBbb=msbm10 
\font\sevenBbb=msbm7 
\font\fiveBbb=msbm5 
\textfont\Bbbfam=\tenBbb 
\scriptfont\Bbbfam=\sevenBbb 
\scriptscriptfont\Bbbfam=\fiveBbb

\newtheorem{theorem}{Theorem}[section] 
\newtheorem{lemma}[theorem]{Lemma}

\theoremstyle{definition}

\newcommand{\s}{\sigma}

\newcommand{\Z}{\mathbb Z}

\def\1{{\mathchoice {1\mskip-4mu\mathrm l}      
{1\mskip-4mu\mathrm l} 
{1\mskip-4.5mu\mathrm l} {1\mskip-5mu\mathrm l}}}

\renewcommand{\qed}{\hfill\ensuremath{\square}}

\renewcommand{\d}{{\rm d}}

\newcommand{\e}   {{\operatorname e }}

\numberwithin{equation}{section}

\begin{document}
\title{Responsive Dormancy of a Spatial Population Among a Moving Trap}
\maketitle
\thispagestyle{empty}
\vspace{-0.5cm}

\centerline{\sc 
Helia Shafigh\footnote{WIAS Berlin, Mohrenstra{\ss}e 39, 10117 Berlin, Germany, {\tt helia.shafigh@wias-berlin.de}} and Leo Tyrpak\footnote{University of Oxford Department of Statistics, 24 St Giles', Oxford OX1 3LB, United Kingdom, {\tt leo.tyrpak@stats.ox.ac.uk}}}
\renewcommand{\thefootnote}{}
\vspace{0.5cm}
\centerline{\textit{WIAS Berlin, and University of Oxford}}

\bigskip


\begin{abstract}
In this paper, we study a spatial model for dormancy in a random environment via a two-type branching random walk in continuous-time, where individuals switch between dormant and active states depending on the current state of a fluctuating environment (responsive switching). The branching mechanism is governed by the same random environment, which is here taken to be a simple symmetric random walk. We will interpret the presence of this random walk as a \emph{trap} which attempts to kill the individuals whenever it meets them. The responsive switching between the active and dormant state is defined so that active individuals become dormant only when a trap is present at their location and remain active otherwise. Conversely, dormant individuals can only wake up once the environment becomes trap-free again.

We quantify the influence of dormancy on population survival by analyzing the long-time asymptotics of the expected population size. The starting point for our mathematical considerations and proofs is the parabolic Anderson model via the Feynman-Kac formula. Specifically, we investigate the quantitative role of dormancy by extending the Parabolic Anderson model to a two-type random walk framework.
 \end{abstract}



\section{Introduction and main results}
\subsection{Biological Motivation} Dormancy is an evolutionary trait that has developed independently across various life forms and is particularly common in microbial communities. To give a definition, we follow \cite{blath} and refer to dormancy as the \emph{ability of individuals to enter a reversible state of minimal metabolic activity}. The collection of all dormant individuals within a population is also often called a \emph{seed-bank}. Maintaining a seed bank reduces the reproduction rate but also decreases resource requirements, making dormancy a viable strategy under unfavourable environmental conditions. Initially studied in plants as a survival strategy (cf.\,\cite{cohen}), dormancy is now recognized as a prevalent trait in microbial communities with significant evolutionary, ecological, and pathogenic implications, serving as an efficient strategy to survive challenging environmental conditions, competitive pressure, or antibiotic treatment. However, it is at the same time a costly trait whose maintenance requires energy and
a sophisticated mechanism for switching between active and dormant states. Moreover, the increased survival rate of dormant individuals must be weighed against their low reproductive activity. Despite its costs, dormancy still seems to provide advantages in variable environments. For a recent overview on biological dormancy and seed-banks we refer to \cite{dormancy}.

The existing stochastic models for dormancy can be roughly categorized into two approaches: population genetics models and population dynamics models. The first approach assumes a constant population size and focuses on the genealogical implications of seed banks, whereas the second typically deals with individual-based modeling via branching process theory. Following a brief example in the book \cite{bookhaccou}, a two-type branching process (without migration) in a fluctuating random environment has been introduced in \cite{blath}, which served as a motivation for this paper. In \cite{blath}, the authors discuss three switching strategies between dormant and active types: \emph{stochastic (or: spontaneous; simultaneous) switching}, \emph{responsive switching} and \emph{anticipatory switching}. In the latter two strategies, individuals adapt to environmental fluctuations by choosing their state (dormant or active) based on environmental cues—for example, with increased reproduction during favourable phases and a larger seed bank during unfavourable ones in the responsive strategy, and vice versa in the anticipatory strategy. The stochastic strategy, in contrast, remains unaffected by environmental changes.

\subsection{Modelling Approach and Goals}  
The aim of this paper is to investigate the responsive switching strategy in order to quantitatively compare the long-term behaviour of populations with and without this dormancy mechanism, in the case where the underlying environment is random and consists of a single moving particle.

Inspired by the Galton-Watson process with dormancy introduced in \cite{blath}, a spatial model for dormancy in a random environment was recently proposed in \cite{shafigh}, where the effect of stochastic dormancy on population growth and survival on $\mathbb{Z}^d$ is quantified via the large-time asymptotics of the expected population size. The random environment driving the population dynamics was modelled using three types of particle systems: a \emph{Bernoulli field of immobile particles}, a \emph{single moving particle}, and a \emph{Poisson field of moving particles}. This framework was further extended in \cite{shafigh25} to incorporate the \emph{simple symmetric exclusion process} as the underlying random environment. To the best of our knowledge, other spatial models of dormancy in random environments within the setting of \emph{population size} models are still absent from the literature. In particular, the \emph{responsive} dormancy strategy has not yet been addressed in the context of spatially structured populations.
 
\subsection{Description of the Model}
In our model, the population resides on $\mathbb{Z}^d$ and consists of two different types $i \in \{0,1\}$ of particles, where we refer to $0$ as \emph{dormant} and to $1$ as \emph{active}. Given random rates $\xi^+(x,t), \xi^-(x,t) \geq 0$, which depend on space and time, we define $\eta(x,i,t)$ to be the number of particles at spatial point $x \in \mathbb{Z}^d$ and state $i$ at time $t \geq 0$, which evolves over time according to the following rules:
\begin{itemize}
\item at time $t=0$, there is only one active particle in $0\in\mathbb{Z}^d$ and all other sites are vacant;
\item all particles act independently of each other;
\item active particles become dormant at rate $s_1(\xi^+-\xi^-)\geq 0$ and dormant particles become active at rate $s_0(\xi^+-\xi^-)\geq 0$;
\item active particles split into two at rate $\xi^+$ and die at rate $\xi^-\geq 0$;
\item active particles jump to one of the neighbour sites with equal rate $\kappa\geq 0$;
\item dormant particles do not participate in branching, dying or migration.
\end{itemize}
By assumption, the initial condition is given by $\eta(x,i,0)=\delta_{(0,1)}(x,i)$. Let us define $\eta(t):=\left\{\eta(x,i,t)\mid (x,i)\in\mathbb{Z}^d\right\}$ as configurations on $\mathbb{N}^{\mathbb{Z}^d\times\{0,1\}}$, representing the number of particles in each point $x\in\mathbb{Z}^d$ and state $i\in\{0,1\}$ at time $t$. Then $\eta=(\eta(t))_{t\geq 0}$ is a Markov process on $\mathbb{N}^{\mathbb{Z}^d\times\{0,1\}}$. However, as we will see later, we will use other methods throughout the paper to describe our population, such that a further formalization of $\eta$ shall not be necessary. In the following, we abbreviate $\xi(x,t):=\xi^+(x,t)-\xi^-(x,t)$ for the \emph{balance} between branching and dying and refer to $\xi$ as the underlying \emph{random environment}. In the following, if we fix a realization of $\xi$, then we will denote by
\begin{align}\label{udef}
u(x,i,t):=u_{\xi}(x,i,t):=\mathbb{E}[\eta(x,i,t)\mid \xi]
\end{align}
the expected number of particles in $x\in\mathbb{Z}^d$ and state $i\in\{0,1\}$ at time $t$ with initial condition 
\begin{align*}
u(x,i,0)=\delta_{(0,1)}(x,i),
\end{align*}
where the expectation is only taken over switching, branching and dying (i.\,e.\, over the evolution of $\eta$ for fixed $\xi$) and not over the random environment $\xi$. If we average over the environment $\xi$ as well, we use the following notation:
\begin{align*}
\left<u(x,i,t)\right>
\end{align*}
as the \emph{annealed} number of particles in $x\in\mathbb{Z}^d$ and in state $i\in\{0,1\}$ at time $t$. 

\subsection{Choice of the Random Environment}
Although one could, in principle, allow $\xi$ to take both negative and positive values—corresponding to death and branching of individuals, respectively—in this paper we restrict ourselves to a \emph{trapping} random environment defined as follows. $Y=(Y(t))_{t\geq 0}$ is a continuous-time simple symmetric random walk on $\mathbb{Z}^d$ that starts in the origin and jumps at total rate $2d\rho$ for some $\rho>0$. We then define $\xi$ as a Markov process on $\{-\gamma,0\}^{\mathbb{Z}^d}$ given by
\begin{align*}
\xi(x,t)=-\gamma\delta_{Y(t)}(x), \qquad x\in\mathbb{Z}^d, t>0,
\end{align*}
where $\gamma>0$ is fixed. In words, the environment contains a single moving particle; whenever it occupies the same site as an individual, that individual is trapped and dies
at rate $\gamma>0$.

\subsection{Main Result}
Recall the number of particles $u(x,i,t)$ in point $x\in\mathbb{Z}^d$ and state $i\in\{0,1\}$ at time $t$, as defined in \eqref{udef}. The quantity we are interested in at most in the current paper is the \emph{annealed expected number of all particles}
\begin{align}\label{grossu}
\left<U(t)\right>:=\sum_{x\in\mathbb{Z}^d}\sum_{i\in\{0,1\}}\left<u(x,i,t)\right>,
\end{align}
which turns into the \emph{annealed survival probability} up to time $t$, if we start with one single particle. Our main result reads:

\begin{theorem}\label{mainresult}
For all $\gamma\in(0,\infty)$ we have that
\begin{align}\label{survival}
\left<U(t)\right>=\left\{\begin{array}{ll}\displaystyle \frac{2(\sqrt{\kappa+\rho}+s_1C_{1,\rho,\kappa,s_1})}{\gamma\sqrt{\pi t}}(1+o(1)), &d=1,\\[13pt]\displaystyle \frac{4(\kappa+\rho)\pi+s_1C_2}{\gamma\log(t)}(1+o(1)), &d=2,\\[13pt]\displaystyle 1-\frac{\gamma \left(G_d(0)+\frac{s_1}{2d(\kappa+\rho)}\right)}{\kappa+\rho +\gamma \left(G_d(0)+\frac{s_1K_d}{2d(\kappa+\rho)}\right)}, &d\geq 3\end{array}\right.
\end{align}
as $t\to\infty$, where $G_d(0)$ is the Green's function of a random walk with total jump rate $2d$, $C_2:=C_{2,s_1,\rho}>0$ and  $K_d:=K_{d,s_1,\rho}>0$ are some positive constants depending on $\rho$ and $s_1$, and
\begin{align*}
    C_{1,\rho,\kappa,s_1}=\frac{1}{\sqrt{\kappa+\rho}\left(\sqrt{s_1^2+4\rho s_1}-s_1\right)}.
\end{align*}
\end{theorem}
We note that the parameter $s_0$ doesn't appear in the main result, which at first sight might seem surprising. We will comment more on this in \ref{section:discussion}.

\subsection{Related results}
The parabolic Anderson model in the absence of switching has attracted considerable interest in recent years and has been extensively studied for various random environments comprising interacting particle systems. For a comprehensive overview of recent developments concerning the parabolic Anderson model, we refer the reader to \cite{PAM}. In this section, we restrict our attention to those results concerning the parabolic Anderson model on $\mathbb{Z}^d$ that are most pertinent to and closely aligned with our framework. The setting of a random walk in the presence of a \emph{single mobile trap} has been investigated in \cite{sw}. In this context, the random walk and the trap independently jump to neighboring lattice sites at rates $\kappa$ and $\rho$, respectively, and the walk is killed at rate $\gamma$ whenever both particles occupy the same site. It was shown in \cite{sw} that the survival probability of the random walk exhibits the following asymptotic behavior:
\begin{align}\label{schnitzler}
\left\{\begin{array}{ll}\displaystyle\frac{2\sqrt{\rho+\kappa}}{\gamma\sqrt{\pi t}}(1+o(1)), &d=1,\\[12pt]\displaystyle\frac{4\pi(\rho+\kappa)}{\gamma\log(t)}(1+o(1)), &d=2,\\[12pt]\displaystyle 1-\frac{\gamma G_d(0)}{\rho+\kappa+\gamma G_d(0)}, &d\geq 3,\end{array}\right.
\end{align}
for $t\to\infty$ where $G_d(0)$ denotes the Green's function of a random walk with jump rate $2d$ evaluated at the origin. Accordingly, in dimensions $d \in \{1,2\}$, the survival probability decays in time, with a constant governed by the parameters $\rho$, $\kappa$, and $\gamma$. In contrast, for $d \geq 3$, the survival probability converges to a strictly positive limit in $(0,1)$, which additionally depends on the expected total time the random walk spends at the origin, as encoded in the Green's function.

More recently, the parabolic Anderson model incorporating a \emph{stochastic dormancy strategy}—in which switching rates between active and dormant states are constant and independent of the environment—has been studied in \cite{shafigh} and \cite{shafigh25} for various specific realizations of the random environment $\xi$. Specifically, the annealed particle density $\langle U(t) \rangle$ has been analyzed in the cases where $\xi$ is given by (1) a Bernoulli field of static particles, (2) a single moving particle, (3) a Poisson field of independently moving particles, and (4) the simple symmetric exclusion process. In the second case—featuring a single mobile trap, which is most directly related to our setting—it has been demonstrated in \cite[Theorem 1.2]{shafigh} that the annealed survival probability decays to zero for $d \in {1,2}$ as $t \to \infty$, and satisfies the asymptotic relation
\medskip
\begin{align}\label{asy2-}
\left\{\begin{array}{ll}\displaystyle\frac{2\sqrt{(s_0+s_1)(s_0(\rho+\kappa)+s_1\rho)}}{s_0\gamma\sqrt{\pi t}}(1+o(1)), &d=1,\\[12pt]\displaystyle\frac{4\pi(\frac{s_1}{s_0}\rho+\rho+\kappa))}{\gamma\log(t)}(1+o(1)), &d=2, \\[12pt]\displaystyle1-\frac{\gamma G_d(0)}{\frac{s_0}{s_0+s_1}\left(\rho+\frac{s_0}{s_0+s_1}\kappa\right)+\gamma G_d(0)}, &d\geq 3\end{array}\right.
\end{align}
as $t\to\infty$, where $G_d$ is the Green's function of a random walk with jump rate $2d$. A comparison of the decay rates and prefactors in \eqref{schnitzler} and \eqref{asy2-} reveals that the incorporation of the stochastic dormancy strategy leads to an enhancement of the survival probability across all spatial dimensions.

\subsection{Relation to the Parabolic Anderson Model and the Feynman-Kac formula}
Consider a one-type branching random walk with exclusively active particles, evolving according to the same dynamics as in our model, except for the switching mechanism, and initiated by a single particle at the origin. It is a well-established result (proved in e.g. \cite{garmol}) that the expected number of particles $u(x,t)$ at spatial point $x$ and time $t$ solves the \emph{Parabolic Anderson model}
\begin{align*}
\left\{\begin{array}{lllr}\frac{\d}{\d t}u(x,t)&=&\kappa\Delta u(x,t) + \xi(x,t)u(x,t), &t>0, x\in\mathbb{Z}^d \\[12pt]u(x,0)&=&\delta_0(x), &x\in\mathbb{Z}^d,
\end{array}\right.
\end{align*}
where $\Delta$ is the discrete Laplacian
\begin{align*}
\Delta f(x):=\sum_{y\in\mathbb{Z}^d, x\sim y}[f(y)-f(x)]
\end{align*}
acting on functions $f:\mathbb{Z}^d\to\mathbb{R}$.
We note that this is a coupled system of ODE's $(u(x,.))_{x\in\Z^d}$ with random coefficients $\xi(x,t)$.
The parabolic Anderson model has been studied intensely during the past years and a comprehensive overview of results can be found in \cite{PAM}. One of the most powerful tools and often the starting point of the analysis of the PAM is the \emph{Feynman-Kac formula}
\begin{align}\label{fkwithout}
u(x,t)=\mathbb{E}_{x}^{X}\left[\exp\left(\int_0^t\xi(X(s),t-s)\,\d s\right)\delta_0(X(t))\right],
\end{align}
where $\mathbb{E}_{x}^{X}$ denotes the expectation with respect to a continuous-time simple symmetric random walk $X$ with start in $x$ and generator $\kappa\Delta$. In other words, the Feynman-Kac formula asserts that the time evolution of all particles can be expressed as an expectation over one single particle moving around according to the same migration kernel and with a \emph{varying mass}, representing the population size. As we can see on the right hand-side of \eqref{fkwithout}, the mass of $X$ changes exponentially depending on the random environment $\xi$. Now, since the Feynman-Kac formula is a powerful tool for the study of the parabolic Anderson model, it is natural to seek an analogous representation for our two-type branching model with switching. To this end, let $Y=(Y(t))_{t\geq 0}$ be a continuous-time simple symmetric random walk on $\mathbb{Z}^d$ with jump rate $2d\rho$ for a constant $\rho>0$ and starting in the origin. For a fixed realization of $Y$, we want to define a joint process $(X,\alpha)$ on $\mathbb{Z}^d\times\{0,1\}$ with the following dynamics: Whenever $\alpha$ is $1$, the process $X$ performs a simple symmetric random walk with total jump rate $2d\kappa$, and it stays still if $\alpha$ is $0$. On the other hand, the dynamics of $\alpha$ is prescribed by the trap $Y$ as well as the walk $X$, in the sense that $\alpha$ may jump from $1$ to $0$ with rate $s_1$, whenever the trap $Y$ and the walk $X$ meet, and jump from $0$ to $1$ with rate $s_0$, whenever $X$ is away from the trap. Further, recall the quantity $u(x,i,t)$ defined in \eqref{udef}, which represents the number of individuals of the population in spatial position $x\in\mathbb{Z}^d$, state $i\in\{0,1\}$ at time $t\geq 0$. Then, given a fixed realization of $Y$, the function $u:\mathbb{Z}^d\times\{0,1\}\times[0,\infty)\to\mathbb{R}$ can be interpreted as the formal solution of the partial differential equation
\medskip
\begin{align}\label{pamswitching}
\left\{\begin{array}{lllr}\frac{\d}{\d t}u(x,i,t)&=&i\kappa\Delta u(x,i,t) + Q u(x,i,t)-i\gamma\delta_{Y(t)}(x)u(x,i,t), &t>0, \\[12pt]u(x,i,0)&=&\delta_{(0,1)}(x,i).
\end{array}\right.
\end{align}
Here, the operator $Q$ is defined as
\begin{align*}
Q u(x,i,t):= s_i(x-Y(t))(u(x,1-i,t)-u(x,i,t)),
\end{align*}
and $s_i(z)$ is defined as
\begin{align*}
s_1(z):=\left\{\begin{array}{ll}s_1, & z=0,\\0, &\text{otherwise},\end{array}\right.\qquad s_0(z):=\left\{\begin{array}{ll}0, & z=0,\\s_0, &\text{otherwise},\end{array}\right.
\end{align*}
for constant rates $s_0,s_1\geq 0$, and
\begin{align*}
\Delta u(x,i,t):=\sum_{y\in\mathbb{Z}^d, x\sim y}[u(y,i,t)-u(x,i,t)]
\end{align*}
(cf.\cite{baran}). Clearly, without $Y$ the process $(X,\alpha)$ is not Markovian. However, as long as we are only interested in the annealed quantity $\left<u(x,i,t)\right>$ after averaging over $Y$ as well, it is sufficient to find a proper formulation for the dynamics of the triple $(Y,X,\alpha)$. To this end, set $Z:=X-Y$. Then we may describe the Markovian dynamics of $(Y,X,\alpha)$ by defining the Markov process $(Z,\alpha)$ which has the generator
\begin{align*}
\bar{\mathcal{L}}f(z,i)=\sum_{y\sim z}(i\kappa+\rho)(f(y,i)-f(z,i))+s_i(z)(f(z,1-i)-f(z,i))
\end{align*} 
for $z\in\mathbb{Z}^d$, $i,j\in\{0,1\}$ and a test function $f:\mathbb{Z}^d\times \{0,1\}\to\mathbb{R}$. Then, we call $(Z,\alpha)$ a \emph{regime-switching random walk} (cf.\,\cite{switching} for the continuous-space version) and the corresponding \emph{Feynman-Kac formula} reads
\begin{align*}
\left<u(x,i,t)\right>=\mathbb{E}_{(x,i)}^{(Z,\alpha)}\left[\exp\left(-\int_0^t\gamma\alpha(s)\delta_{0}(Z(s))\,\d s\right)\delta_{(0,1)}(Z(t),\alpha(t))\right],
\end{align*}
where $\mathbb{E}_{(x,i)}^{(Z,\alpha)}$ denotes the expectation over $(Z,\alpha)$ starting from $(x,i)$ (cf.\,\cite{baran}). To obtain the total number $\left<U(t)\right>$ of particles at time $t$, we can sum over all states $(x,i)$ and use a time-reversal, yielding
\begin{align}\label{Uresponsive}
\left<U(t)\right>=\mathbb{E}_{(0,1)}^{(Z,\alpha)}\left[\exp\left(-\gamma\int_0^t\delta_{(0,1)}(Z(s),\alpha(s))\,\d s\right)\right],
\end{align}
where $\mathbb{E}_{(x,i)}^{(Z,\alpha)}$ denotes the expectation with respect to the Markov process $(Z,\alpha)$ started in $(x,i)$. Thus, the study of the first moment of our two-type branching process can be reduced to the analysis of only one particle with the same migration, branching and switching rates. Note that, as we will only consider \emph{traps} in this paper, the quantity \eqref{Uresponsive} lies in $[0,1]$ in this setting and represents the annealed survival probability of $X$ up to time $t$.

\subsection{Discussion}\label{section:discussion}
In this section, we analyse how the incorporation of a responsive dormancy strategy, as introduced in our model, influences the long-time dynamics of the population. To begin, a direct comparison between the known asymptotics \eqref{schnitzler} and our result \eqref{survival} reveals that, across all spatial dimensions, the inclusion of responsive dormancy increases the survival probability relative to corresponding models without dormancy. In the one-dimensional setting, the improvement due to dormancy is quantified by the additive term
\begin{align}\label{add1}
\frac{s_1}{\gamma\sqrt{(\kappa\rho)\pi t}\left(\sqrt{s_1^2+4\rho s_1}-s_1\right)},
\end{align}
which vanishes in the limit $s_1\to 0$, thereby recovering the known asymptotics \eqref{schnitzler}. A straightforward calculation shows that this term is strictly increasing in $s_1$, confirming that the responsive strategy becomes increasingly beneficial with stronger dormancy. The more subtle question, however, is how responsive dormancy compares to stochastic dormancy. While both strategies yield identical asymptotic survival probabilities as $s_1\to 0$ (converging to \eqref{schnitzler}), their large-$s_1$ behaviour is different. In the stochastic case, the leading-order prefactor becomes $\frac{2\sqrt{\rho}s_1}{s_0}$, whereas for the responsive strategy we obtain $\frac{s_1}{\rho\sqrt{\kappa+\rho}}$. Comparing these expressions shows that the responsive strategy yields a higher survival probability for large $s_1$ if and only if 
\begin{align*}
    s_0>2\rho^{\frac{3}{2}}\sqrt{\kappa+\rho}.
\end{align*}
We do not know of any heuristic explanation for this result.
Thus, in one dimension, there is a threshold behaviour governed by the interplay between the reactivation rate $s_0$, the dormancy intensity $s_1$ and the motion parameters $\kappa, \rho$: the responsive strategy outperforms the stochastic one when reactivation is sufficiently frequent relative to the underlying mobility dynamics. A similar effect is also observed in two dimensions. Here, the survival probability under the responsive strategy is enhanced by an additive term of the form $\frac{s_1C_2}{\gamma\log(t)}$, where $C_2$ depends on $s_1$ and $\rho$ only. Although this additive term is always positive for $s_1>0$ and vanishes for $s_1=0$, yielding the known asymptotics \eqref{schnitzler} in this case, it is not immediately clear whether the survival probability is monotonic in. To compare the responsive and stochastic strategies, we have to compare the additive terms $\frac{4\pi s_1\rho}{s_0}$ and $s_1 C_2$ respectively, as seen from \eqref{asy2-} and \eqref{mainresult} respectively, which yields a simple criterion: the responsive strategy provides a higher survival probability if and only if 
 \begin{align*}
     C_2>\frac{4\pi\rho}{s_0}.
 \end{align*}
 Thus, the relative effectiveness of the two strategies depends on the magnitude of $C_2$ appearing in \eqref{survival}, which is unfortunately not explicit.
 However some intuition is that the survival probability in the responsive dormancy case does not depend on $s_0$ whereas the survival probability of stochastic strategy is worse when $s_0$ is very large as that means individuals leave dormancy very fast.
 For higher dimensions, a direct comparison of the asymptotic survival probabilities from \eqref{survival} and \eqref{schnitzler} leads to the inequality
\begin{align}\label{3diff}
\frac{G_d(0) + \frac{s_1}{2d(\kappa + \rho)}}{\kappa + \rho + \gamma \left(G_d(0) + \frac{s_1 K_d}{2d(\kappa + \rho)}\right)} < \frac{G_d(0)}{\kappa + \rho + \gamma G_d(0)},
\end{align}
which holds for all $s_1>0$. This confirms that the responsive strategy always improves survival probability when dormancy is present. Furthermore, equality is achieved at $s_1=0$, as expected. The derivative of the left-hand side of \eqref{3diff} with respect to $s_1$ is negative under the condition 
\begin{align*}
\kappa+\rho<\gamma G_d(0)(K_d-1)    
\end{align*}
implying that, in this parameter regime, the survival probability is strictly increasing in $s_1$. In addition, while both strategies yield identical long-time survival probabilities in the limit $s_1\to 0$, their behaviour as $s_1\to\infty$ differs: While the stochastic survival probability tends to zero, the responsive one converges to
\begin{align*}
    1-\frac{1}{K_d}.
\end{align*}
Hence, the responsive dormancy strategy consistently outperforms the stochastic one for large dormancy intensities $s_1$.

One interesting—and perhaps surprising—difference between our result and previous ones, both without and with stochastic dormancy, is that, as seen from \eqref{mainresult}, the asymptotic survival probability in our model appears to be entirely independent of the reactivation rate $s_0$, whereas $s_0$ features explicitly in all asymptotic expressions derived in \cite{shafigh}. This phenomenon can be interpreted as follows: The key difference between the present model and earlier results such as those in equation \eqref{asy2-} lies in how the switching mechanism between active and dormant states interacts with the environment. In previous work, the switching rates $s_0$ and $s_1$ are constant and independent of the particle's position relative to the trap. As a result, particles may become active even while sitting directly on the trap, exposing themselves to an immediate risk of death. In this setting, both $s_0$ and $s_1$ influence the overall survival probability, since they determine the average fraction of time a particle spends in the active state — the only state in which it can both reproduce and die. Consequently, the balance between activation and dormancy directly affects both the potential for population growth and the likelihood of extinction, making both switching rates relevant to the asymptotics of the survival probability. In contrast, the present model introduces a structured, environment-dependent switching mechanism: particles become dormant only when located at the trap, and they can reactivate only when they are away from it. This spatial constraint introduces a protective behaviour — particles automatically hide when in danger and only reemerge in safe regions. As a result, the reactivation rate $s_0$ has no direct impact on the likelihood of survival, since activation never occurs in risky locations. Instead, only the dormancy rate $s_1$, which governs how efficiently particles avoid danger by entering the dormant state, plays a role in the long-term asymptotics. Since dormant particles can not be killed, and reactivation happens only away from the trap, the extinction dynamics are entirely driven by the behaviour of the active particles, and the protective mechanism encoded in $s_1$ determines the survival probability — leaving $s_0$ absent from the final results.

\subsection{Outline}
The remainder of the paper is structured as follows. In Section 2, we derive a representation of the survival probability based on the number of times the random walk is trapped while the particle is in its active state. This representation is further related to the distribution of the waiting times between successive trapping events, for which we establish a precise asymptotic characterization. Section 3 is devoted to the proof of Theorem 1.1.

\section{Preparatory Results}
Before we proceed with the proof of Theorem \ref{mainresult}, we first need some preparations:
\subsection{Regeneration times}
As the definition of the responsive model suggests, the dynamics of the process $(Z,\alpha)$ are more complicated than in the case of the stochastic dormancy, so that our proof techniques from \cite{shafigh} are not applicable. Hence, we need to find a new way to analyze the asymptotics of the survival probability \eqref{Uresponsive} for $t\to\infty$. To this end, we will make use of renewal theory by noting that, starting from $(0,1)$, every time the process $(Z,\alpha)$ comes back to $(0,1)$, it \emph{regenerates} in the sense that it forgets its past. Define the sequence of stopping times $(\tau_n)_{n\in\mathbb{N}}$ by $\tau_0:=0$ and
\begin{align*}
\sigma_n:=&\inf\left\{t\geq \tau_{n-1}:(Z(t),\alpha(t))\neq(0,1)\right\},\\
\tau_n:=&\inf\left\{t\geq \sigma_n:(Z(t),\alpha(t))=(0,1)\right\},\qquad n\in\mathbb{N}, n\geq 1. 
\end{align*}
Note: Define the time $\sigma_n$ when we leave $(0,1)$.\\
If we let $Z_n=\tau_n-\tau_{n-1}$, then random variable $Z_1,...,Z_n$ are independent and identically distributed with some distribution function $F$, the \emph{inter-arrival distribution}, and satisfy $\tau_n=\sum_{i=1}^nZ_i$. Define the counting process $(N_t)_{t\geq0}$ by
\begin{equation}
    N_t=|\{n\geq1:\tau_n\leq t\}|, \qquad t>0,
\end{equation}
Then
\begin{align*}
    \left<U(t)\right>=\mathbb{E}\left[\exp\left(-\gamma\sum_{i=1}^{N_t}Y_i\right)\right],
\end{align*}
where $Y_i$ is the time spent at the state $(0,1)$ at the $i$-th visit. We note that $Y_i$, $i\geq 1$, are independent and exponentially distributed with parameter $2d(\kappa+\rho)+s_1$.
Therefore,
\begin{align}\label{defGmu}
    \left<U(t)\right>=\sum_{n\geq 1}\mathbb{P}[N_t=n]\mu^n=:G_\mu(t)
\end{align}
for $\mu:=\frac{2d(\kappa+\rho)+s_1}{2d(\kappa+\rho)+s_1+\gamma}$.
Hence, the study of the survival probability reduces to the study of the counting process $(N_t)_{t\geq 0}$. As a first step, we clarify the relation between the counting process $(N_t)_{t\geq 0}$ and the inter-arrival times $Z_n$, $n\in\mathbb{N}$, which exhibit dimension-dependent behaviour. Our first lemma states the desired relation in the recurrent dimensions:

\begin{lemma}\label{G}
Denote by
\begin{align*}
\overline{g}(\lambda)=\int_0^{\infty}\e^{-\lambda t}\d g(t)
\end{align*}
the Laplace-Stieltjes transform of a real-valued function $g$ and by
\begin{align*}
\widehat{g}(\lambda)=\int_0^{\infty}\e^{-\lambda t}g(t)\d s
\end{align*}
its Laplace transform. Recall the distribution function $F(t)=\mathbb{P}(Z_1\leq t)$ of $Z_1$, as well as $G_\mu$ defined in \eqref{defGmu}. Then, for all $\lambda>0$,
\begin{align}\label{survivallaplace}
    \widehat{G}_\mu(\lambda)=\frac{\mu}{\lambda}\cdot\frac{1-\lambda\overline{F}(\lambda)}{1-\mu\lambda\overline{F}(\lambda)}. 
\end{align}
\end{lemma}
\begin{proof}[Proof]
    Conditioned on $Z_1$, we have
    \begin{align*} 
    G_\mu(t)=\sum_{n\geq 1}\mu^n\left(\mathbb{P}(N_t=n\mid Z_1>t)\mathbb{P}(Z_1>t)+\mathbb{P}(N_t=n\mid Z_1\leq t)\mathbb{P}(Z_1\leq t)\right).
    \end{align*}
On the event $\{Z_1>t\}$ we have
\begin{align*}
    \sum_{n\geq 1}\mu^n\mathbb{P}(N_t=n\mid Z_1>t)=\mu\mathbb{P}(N_t=1\mid Z_1>t)=\mu,
\end{align*}
as we then hit $(0,1)$ exactly once, namely in the starting point. However, on $\{Z_1\leq t\}$ we restart the counting process at $Z_1$ and so, given $\{Z_1\leq t\}$, $N_t$ is equal to $1+N_{t-Z_1}$ in distribution. Thus, 
\begin{align*}
    \sum_{n\geq 1}\mu^n\mathbb{P}(N_t=n\mid Z_1\leq t)=\mu\int_0^t\sum_{n\geq 1}\mu^n\mathbb{P}(N_{t-s}=n)\d F(s).
\end{align*}
This leads to the recursion formula
\begin{equation*}
    G_\mu(t)=\mu(1-F(t))+\mu\int_0^tG_\mu(t-s)\d F(s).
\end{equation*}
Denoting $f(s):=\frac{\d}{\d s}F(s)$ and taking the Laplace transform on both sides, we obtain
\begin{equation*}
    \widehat{G}_\mu(\lambda)=\frac{\mu}{\lambda}-\mu\overline{F}(\lambda)+\mu\widehat{G}_\mu(\lambda)\hat{f}(\lambda).
\end{equation*}
Noting that $\widehat{f}(\lambda)=-F(0)+\lambda\overline{F}(\lambda)=\lambda\overline{F}(\lambda)$, and solving for $\widehat{G}_\mu(\lambda)$, we obtain the claim. 
\end{proof}
\subsection{Distribution of $Z_1$}
As seen from the relation \eqref{survivallaplace}, the distribution function $F$ of the inter-arrival time $Z_1$ is crucial for calculating the long-time asymptotics of the survival probability. Hence, in the next step we aim to study the behaviour of $Z_1$ in more detail. This will be done by decomposing the random time $Z_1$ into a sum of independent random variables. We begin with a heuristic description of the decomposition. Suppose $(Z,\alpha)$ is at state $(0,1)$ and
let $A$ be the event that the walk $Z$ jumps away from $0$, say to $e_1$ without loss of generality, before transitioning to the dormant state or being killed. The probability of this event is given by $\frac{2d(\kappa+\rho)}{2d(\kappa+\rho)+s_1}$. 
On the event $A$, the walk $Z$ behaves like a simple symmetric random walk with total jump rate $2d(\kappa+\rho)$, so that $Z_1$ corresponds to the first hitting time of the origin by this walk, when starting in $e_1$. We will denote this first hitting time by $R^{2d(\kappa+\rho)}_{e_1,0}$ in the following. On the complementary event $A^c$, the walk becomes dormant before making a jump. In this case, the only possible transition is that $Z$ moves away from $0$, say to $e_1$ w.l.o.g, according to a simple symmetric random walk with jump total rate $2d\rho$. The process may then either return to the origin before becoming active again, initiating another excursion, or it may be reactivated at a non-zero site. Consequently, there is a geometrically distributed number $N$ of independent excursions from $0$ to $0$ with jump rate $2d\rho$, each conditioned to complete before the exponential waiting time until reactivation expires. In the following, we denote the duration of the $i$-th such excurison by $R^{2d\rho, i}_{0,0}$. Now, assume that the walk wakes up before hitting $0$ again, and denote by $Y$ its location in the moment it gets active, which is a random variable on $\mathbb{Z}^d$. From now on, $Z$ performs a simple symmetric random walk with jump rate $2d(\kappa+\rho)$ before hitting $0$ again. We denote by $R_{0,Y}^{2d(\kappa+\rho)}$ the random variable representing the duration of this return time. This leads to the decomposition
\begin{align}\label{decomposition}
    Z_1\overset{d}= &\left(E^{2d(\kappa+\rho)+s_1}+R_{e_1,0}^{2d(\kappa+\rho)}\right)1_A+\\
    &\left(\sum_{i=1}^N\left(E^{2d\rho+s_1,i}+R_{e_1,0}^{2d\rho,i}\right)+E^{s_1}+R_{Y,0}^{2d(\kappa+\rho)}\right)1_{A^c},
\end{align}
where we wrote $E^\beta$ for an exponential-distributed random variable with parameter $\beta$. 
Some notes about this decomposition.
\begin{itemize}
    \item Here $R_{e_1,0}^{2d\rho,i}$ are conditioned to be smaller than $E^{s_1}$.
    \item The exponential random variable that appears just after that $E^{s_1}$ is conditioned to be smaller than an independent copy of $R_{e_1,0}^{2d\rho}$, however the duration of that exponential random variable will not be important.
    \item This means that $Y$ is the distribution of a continuous time random walk after $E^{s_1}$ time conditional on $E^{s_1}<R_{e_1,0}^{2d\rho}$, but this conditioning is equivalent to conditioning $Y$ on not hitting $0$.
\end{itemize}
The advantage of this representation is that, since we are summing independent random variables, it will be easier to extract the Laplace transform. In the following, we will denote by
\begin{align}
    f(x,\lambda):=\mathbb{E}\left[e^{-\lambda R^{2d(\kappa+\rho)}_{x,0}}\right]
\end{align}
the Laplace transform of the hitting time of $0$ by a simple symmetric random walk with total jump rate $2d(\kappa+\rho)$ starting in $x\in\mathbb{Z}^d$, which exhibits dimension-dependent behaviour. To simplify notation, we suppress the dependence on $d$ in the abbreviation $f(x,\lambda)$. The following lemma quantifies the asymptotic behaviour of the Laplace transform of $Z_1$ in terms of the function $f$:
\begin{lemma}\label{Zrec}
    Let $\mathbb{E}^Y$ denote the expectation with respect to the random variable $Y$. Then
    \begin{align*}
        \mathbb{E}\left[\exp(-\lambda Z_1)\right]=\frac{2d(\kappa+\rho)}{s_1+2d(\kappa+\rho)}f(e_1,\lambda)+\frac{s_1}{s_1+2d(\kappa+\rho)}\mathbb{E}^Y\left[f(Y,\lambda)\right]+O(\lambda)
    \end{align*}
    as $\lambda\to 0$.
\end{lemma}
\begin{proof}[Proof]
   The proof relies on the representation \eqref{decomposition}. Taking the Laplace transform on both sides, we obtain for the part including the event $A$ that
   \begin{align*}
       \mathbb{E}\left[\e^{-\lambda \left(E^{2d(\kappa+\rho)+s_1}+R_{e_1,0}^{2d(\kappa+\rho)}\right)}1_A\right]&=\frac{2d(\kappa+\rho)}{2d(\kappa+\rho)+s_1}\frac{2d(\kappa+\rho)+s_1}{2d(\kappa+\rho)+s_1+\lambda}f_{e_1}(\lambda)
       \\&=\frac{2d(\kappa+\rho)}{2d(\kappa+\rho)+s_1}f_{e_1}(\lambda)+O(\lambda)
   \end{align*}
as $\lambda\to 0$, since, as we will see in Lemma \ref{d=1} and Lemma \ref{d=2} for dimension $1$ and $2$ respectively, the highest order terms in dimension $1$ and $2$ are $\sqrt{\lambda}$ and $\frac{1}{\log(1/\lambda)}$ respectively, so that terms of order $O(\lambda)$ are negligible. Regarding the term on the event $A^c$ appearing on the right hand-side of \eqref{decomposition}, we have
   \begin{align*}
       \mathbb{E}&\left[\e^{-\lambda \left(\sum_{i=1}^NR_{0,0}^{2d\rho,i}+E^{s_1}+R_{Y,0}^{2d(\kappa+\rho)}\right)}1_{A^c}\right]
       \\&=\frac{s_1}{2d(\kappa+\rho)+s_1}\mathbb{E}\left[\e^{-\lambda\sum_{i=1}^NR_{0,0}^{2d\rho,i}}\right]\mathbb{E}\left[\e^{-\lambda E^{s_1}}\right]\mathbb{E}\left[\e^{-\lambda R_{Y,0}^{2d(\kappa+\rho)}}\right],
   \end{align*}
   where we note that $\mathbb{E}\left[\e^{-\lambda E^{s_1}}\right]=\frac{s_1}{s_1+\lambda}=1+O(\lambda)$, $\lambda\to 0$, so will not contribute to the Laplace transform for small $\lambda$. Moreover, note that for a geometric-distributed random variable $G$ with density $(1-p)^kp$ for $k\geq0$ and independent and identically distributed random variables $X_i$ we have
\begin{align*}
    \mathbb{E}\left[\exp{(-\lambda\sum_{i=1}^GX_i)}\right]=\frac{p}{1-(1-p)\mathbb{E}[\exp{(-\lambda X_1)}]}.
\end{align*}
Moreover,
\begin{align*}
    \mathbb{E}[\exp(-\lambda R_{0,0}^{2d\rho})\mid R_{0,0}^{2d\rho}\leq  E^{s_0})]&\geq \mathbb{E}[\exp(-\lambda E^{s_0})]=\frac{s_0}{s_0+\lambda}
    \\&=1-\frac{\lambda}{s_0}+o(\lambda)=1+O(\lambda)
\end{align*}
as $\lambda\to 0$, which together with the boundedness of the expectation by $1$ gives 
\begin{align*}
\mathbb{E}[\exp(-\lambda R_{0,0}^{2d\rho})\mid R_{0,0}^{2d\rho}\leq  E^{s_0})]=1+O(\lambda), \qquad\lambda\to 0.
\end{align*}
Hence, conditioned on $B:=\{R_{0,0}^{2d\rho,i}\leq E^{s_0} \text{ for all } i\leq N \}$,
   \begin{align*}
       \mathbb{E}\left[\e^{-\lambda\sum_{i=1}^NR_{0,0}^{2d\rho,i}}\mid B \right]&=\frac{\mathbb{P}(R_{0,0}^{2d\rho}> E^{s_0})}{1-\mathbb{P}(R_{0,0}^{2d\rho}\leq  E^{s_0})\mathbb{E}[\e^{-\lambda R_{0,0}^{2d\rho}}\mid R_{0,0}^{2d\rho}\leq  E^{s_0}]}
       \\&=\frac{\mathbb{P}(R_{0,0}^{2d\rho}> E^{s_0})}{1-\mathbb{P}(R_{0,0}^{2d\rho}\leq  E^{s_0})(1+O(\lambda))}=1+O(\lambda)
   \end{align*}
as $\lambda\to 0$. Therefore, this term will not contribute at all, which matches the intuition, as we are summing up a geometric number of random variables all conditioned to be smaller than an exponential one. Combining all and noting that 
\begin{align*}
    \mathbb{E}\left[\exp(-\lambda R_{Y,0}^{2d(\kappa+\rho)})\right]=\mathbb{E}^Y[f(Y,\lambda)],
\end{align*}
where $\mathbb{E}^Y$ shall emphasise the expectation with respect to $Y$, we obtain the claim. 
\end{proof}
Note, that the asymptotic behaviour of $f(e_1,\lambda)$ is already known in the literature (see e.\,g.\, \cite{lawler}) and satsifies
\begin{align}
    f(e_1,\lambda)\sim\left\{\begin{array}{ll}1-\sqrt{\frac{\lambda}{\kappa+\rho }}, &d=1,\\[12pt]1-\frac{\pi}{\log(1/\lambda)}, &d=2,\end{array}\right.
\end{align}
as $\lambda\to 0$. Due to its dimension dependent nature, distinct proof techniques are required for each case. However, the following result, which holds in all dimensions, will play a key role in our analysis:
\begin{lemma}\label{E[Y]}
Let $G$ be a geometric random variable with parameter $p=\frac{s_1}{s_1+2\rho}$, $\tilde{X}$ a discrete-time simple symmetric random walk starting in $e_1$, and $\tau_0$ the first hitting time of the origin by $\tilde{X}$. Then, for any function $h:\mathbb{Z}^d\to[0,\infty)$ which is harmonic on $\mathbb{Z}^d\setminus\{0\}$,
    \begin{align}
        \mathbb{E}^Y[h(Y)]=\frac{h(e_1)}{\mathbb{P}(G<\tau_0)}.
    \end{align}
\end{lemma}
\begin{proof}
    We first note that
     \begin{align*}
      \mathbb{P}^Y(Y=y)=\mathbb{P}_1^X(X_{E^{s_1}}=y\mid E^{s_1}<\tau_0)= \mathbb{P}_1^{\tilde{X}}(\tilde{X}_{G}=y\mid G<\tau_0).
  \end{align*}
  Indeed, by discretising time we obtain
      \begin{align*}
        \mathbb{P}_1^X(X_{E^{s_1}}=y\mid E^{s_1}<\tau_0)&=s_1\int_0^t e^{-s_1 t}\sum_{n=0}^{\infty}\mathbb{P}_1^{\tilde{X}}(\tilde{X}_n=y\mid \tau_0>n)\frac{(2\rho t)^n}{n!}e^{-{2\rho t}}\d t
        \\&=s_1\sum_{n=0}^{\infty}\mathbb{P}_1^{\tilde{X}}(\tilde{X}_n=y\mid \tau_0>n)\frac{(2\rho t)^n}{n!}\cdot\frac{n!}{(2\rho +s_1)^{n+1}}
        \\&=\sum_{n=0}^{\infty}\mathbb{P}_1^{\tilde{X}}(\tilde{X}_n=y\mid \tau_0>n)\left(\frac{2\rho}{2\rho+s_1}\right)^n\frac{s_1}{2\rho+s_1}
        \\&=\mathbb{P}_1^{\tilde{X}}(\tilde{X}_{G}=y\mid G<\tau_0).
      \end{align*}
    Therefore, we will consider $\tilde{X}$ from now on instead of the continuous-time version, and write $\mathbb{P}_{e_1}:=\mathbb{P}_{e_1}^{\tilde{X}}$ for its distribution starting from $e_1$. Now, as $h$ is a harmonic function for $x\neq 0$, $h(\tilde{X}_n)$ is a non-negative martingale conditioned up to the first hitting time $\tau_0$ (on the event $\{n<\tau_0\}$). More precisely, $(h(\tilde{X}_{n\wedge \tau_0}))_n$ is a martingale w.r.t. to the canonical filtration of the random walk $\tilde{X}_n$ and therefore,
  \begin{align*}
      \mathbb{E}_{e_1}[h(\tilde{X}_{n\wedge\tau_0})]=h(e_1), \qquad n\in\mathbb{N}. 
  \end{align*}
  Thus,
  \begin{align*}
      \mathbb{E}[h(\tilde{X}_G)1_{\{G<\tau_0\}}]&=\sum_{n\geq 0}\mathbb{P}(G=n)\mathbb{E}_{e_1}[h(\tilde{X}_n)1_{\{n<\tau_0\}}]
      \\&=h(e_1)\sum_{n\geq 0}\mathbb{P}(G=n)=h(e_1).
  \end{align*} 
Further, note that
  \begin{align*}
      \mathbb{E}[h(Y)]=\mathbb{E}_{e_1}[h(\tilde{X}_G)|G<\tau_0]=\frac{\mathbb{E}_{e_1}[h(\tilde{X}_G)1_{\{G<\tau_0)\}}]}{\mathbb{P}(G<\tau_0)}.
  \end{align*}
\end{proof}

One of the fruits of the previous lemma is the following result concerning the distribution of $Z_1$ in dimension $d=1$:
\begin{lemma}\label{d=1}
    Let $d=1$. Then
    \begin{align*}
        \mathbb{E}[\exp(-\lambda Z_1)]=1-\frac{2\left(\kappa+\rho+s_1C_{1,\rho,s_1}\right)\sqrt{\lambda}}{\sqrt{\kappa+\rho}(s_1+2(\kappa+\rho))}+O(\lambda), \qquad \lambda\to 0,
    \end{align*}
    where $C_{1,\rho,s_1}=\frac{1}{\sqrt{s_1^2+4\rho s_1}-s_1}$.
\end{lemma}
\begin{proof}[Proof]
First, note that if we start in $x$ and hit the origin for the first time at some instant, then this first hitting time is the sum of the first hitting time of $x-1$, and the first hitting time of $x-2$ restarting from $x-1$, as well as the first hitting time of $x-3$ restarting from $x-2$ and so on, so that, due to the independence of all these steps, 
  \begin{align}\label{e_1x}
      R_{x,0}^{2(\kappa+\rho)}\overset{d}=\sum_{i=1}^{|x|}R_i,
  \end{align}
 where $R_i$, $i=1,\cdots,|x|$, are identital copies of $R_{1,0}^{2(\kappa+\rho)}$.
Consequently, $f(x,\lambda)=f(1,\lambda)^{|x|}$. Denoting by $G_Y(s)=\mathbb{E}^Y[s^Y]$ the generating function of $Y$, we obtain
  \begin{align}\label{Ef}
\mathbb{E}^Y[f(Y,\lambda)]&=\mathbb{E}^Y[f(e_1,\lambda)^{\|Y\|}]=G_Y(f(e_1,\lambda))\\
&=1+(f(e_1,\lambda)-1)G_Y'(1)+O(\lambda)\nonumber
   \\&=1-\frac{G_Y'(1)}{\sqrt{\kappa+\rho}}\sqrt{\lambda}+O(\lambda),
  \end{align}
  by Taylor approximation of $G_Y$ around $1$. Thus, we have to calculate $G_Y'(1)=\mathbb{E}[Y]$, which was defined as the position of a random walk with jump rate $2\rho$, starting from $1$, after an exponential time $E^{s_1}$ with parameter $s_1$, and conditioned on not hitting $0$ up to that time. Applying Lemma \ref{E[Y]} to the harmonic function $h:\mathbb{Z}^+\to\mathbb{Z}^+$, $h(z)=z$, yields
  \begin{align*}
      \mathbb{E}[Y]=\frac{1}{\mathbb{P}(G<\tau_0)}.
  \end{align*}
  where $G$ is geometric random variable of parameter $\frac{2\rho}{s_1+2\rho}$. Therefore,
\begin{align*}
    \mathbb{P}(G<\tau_0)=1-\mathbb{E}\left[\left(\frac{2\rho}{s_1+2\rho}\right)^{\tau_0}\right]=\frac{\sqrt{s_1^2+4\rho s_1}-s_1}{2\rho},
\end{align*}
where in the last step we used a known identity for generating functions of first hitting times in dimension $d=1$. Thus,
\begin{align*}
    \mathbb{E}[Y]=\frac{2\rho}{\sqrt{s_1^2+4\rho s_1}-s_1}
\end{align*}
and by Lemma \ref{Zrec},
   \begin{align*}
    \mathbb{E}[\exp(-\lambda Z_1)]&=1-\frac{2(\kappa+\rho)}{s_1+2(\kappa+\rho)}\sqrt{\frac{\lambda}{\kappa+\rho}}
    \\&-\frac{s_1}{s_1+2(\kappa+\rho)}\frac{2\rho\sqrt{\lambda}}{\sqrt{\kappa+\rho}\left(\sqrt{s_1^2+4\rho s_1}-s_1\right)}+O(\lambda)
    \\&=1-\sqrt{\lambda}\frac{2\left(\kappa+\rho+s_1\sqrt{\kappa+\rho}C_{1,\rho,\kappa,s_1}\right)}{\sqrt{\kappa+\rho}(s_1+2(\kappa+\rho))}+O(\lambda),
  \end{align*}
  as $\lambda\to 0$.
\end{proof}
Next, we investigate the asymptotics of the return time $Z_1$ in dimension $d=2$:
\begin{lemma}\label{d=2}
    Let $d=2$, $q=\frac{4\rho}{s_1+4\rho}$, $\phi(k):=\frac{1}{2}(\cos(k_1)+\cos(k_2))$, and
\begin{align*}
    G_2(x,\lambda):=\int_{[-\pi,\pi]^2}\frac{e^{ik\cdot x}}{\lambda+1-\phi(k)}\frac{dk}{(2\pi)^2}
\end{align*}
   the Green's generating function, where the subscript $2$ refers to the dimension. Then
    \begin{align*}
        \mathbb{E}[\exp(-\lambda Z_1)]=1-\frac{\pi}{(s_1+4(\kappa+\rho))\log(1/\lambda)}\left(4(\kappa+\rho)+\frac{s_1G_{2}(0,q)}{G_2(0,1)+G_{2}(e_1,q)}\right)+O(\lambda),
    \end{align*}
    as $\lambda\to 0$.
    \end{lemma}
\begin{proof}
         We will use many results and notations from \cite[Chapter 4]{lawlerlimic}. Note that $f(x,\lambda)$ solves the Poisson equation $(\Delta_x-\lambda)f(.,\lambda)\equiv0$ on $\mathbb{Z}^2\setminus\{0\}$ with boundary condition $f(0,\lambda)=1$. Therefore, by uniqueness of bounded solutions to the Poisson equation, it has the explicit representation $f(x,\lambda)=\frac{G_\lambda(x)}{G_\lambda(0)}$. Let
           \begin{align}\label{defa(x)}
      a(x)=\int_{[-\pi,\pi]^2}\frac{1-e^{ik\cdot x}}{1-\phi(k)}\frac{\d k}{(2\pi)^2}. 
    \end{align}
Then $G_2(x,\lambda)=G_2(0,\lambda)-a(x)-E_\lambda(x)$. Thus, $f(x,\lambda)$ can be written as
    \begin{align}\label{eq:fx_lambda}
      f(x,\lambda)= 1-\frac{a(x)}{G_2(0,\lambda)}+\frac{E_\lambda(x)}{G_2(0,\lambda)}.
    \end{align}
As we will prove in the Appendix (Lemma \ref{lem:asymptotics_G(0)}), we have
\begin{align*}
    G_2(0,\lambda)=\frac{1}{\pi}\log(1/\lambda)+O(1), \qquad\lambda\to 0.
\end{align*}
Moreover, by Lemma \ref{lem:error_E_goes0} we know that $\mathbb{E}[E_\lambda(Y)]\to 0$, as $\lambda\to 0$, since $Y$ jumps only a geometric number of times and so satisfies the hypothesis of finite variance. Taking expectation with respect to $Y$ in (\ref{eq:fx_lambda}) and using the expansion of $G(0,\lambda)$ and the convergence of error $E_\lambda(x)$ we obtain
\begin{align}\label{F2}
    \mathbb{E}[f(Y,\lambda)]=1-\frac{\pi\mathbb{E}_1[a(Y)]}{\log(1/\lambda)}+o\left(\frac{1}{\log(1/\lambda)}\right), \qquad\lambda\to 0.
\end{align} 
Next, as $a$ is a harmonic function for $x\neq 0$, Lemma \ref{E[Y]} asserts that
  \begin{align}\label{a2}
      \mathbb{E}[a(Y)]=\frac{a(e_1)}{\mathbb{P}(G<\tau_0)}=\frac{1}{\mathbb{P}(G<\tau_0)},
  \end{align}
  such that it only remains to determine $\mathbb{P}(G<\tau_0)$. Let
  \begin{align}\label{greenres}
    G^{(0)}_{d}(x-y,\lambda):&=\sum_{n\geq 0}\lambda^n\mathbb{P}_x(\tilde{X}_n=y,\tau_0>n)
  \end{align}
  denote the Green's generating function of $\tilde{X}$ conditioned on not having hit $0$. Then
  \begin{align*}
      \mathbb{P}(G<\tau_0)&=p\sum_{n\geq 0}(1-p)^n\mathbb{P}_{e_1}(\tau_0>n)=p\sum_{n\geq 0}\sum_{y\neq 0}(1-p)^n\mathbb{P}_{e_1}(\tilde{X}_n=y,\tau_0>n)
      \\&=p\sum_{y\neq 0}G_{d}^{(0)}(e_1-y, 1-p).
  \end{align*}
  We decompose the generating Green's function $G_d(\cdot,\lambda)$ into paths avoiding zero and paths hitting zero. For the latter, we can split the first visit to $0$ at some time $k$ and restart from the $0$ to obtain
 \begin{align*}
   G_d(x-y,\lambda)&=G_d^{(0)}(x-y,\lambda)+\sum_{k\geq 1}\lambda^k\mathbb{P}_x(\tau_0=k)G_d(y,\lambda)
   \\&=G_d^{(0)}(x-y,\lambda)+G_d(y,\lambda)\mathbb{E}_x[\lambda^{\tau_0}1_{\{\tau_0<\infty\}}]
\\&=G_d^{(0)}(x-y,\lambda)+G_d(y,\lambda)\frac{G_d(x,\lambda)}{G_d(0,\lambda)}
 \end{align*}
  and hence
  \begin{align*}
   G_d^{(0)}(x-y,\lambda)=G_d(x-y,\lambda)-\frac{G_d(x,\lambda)G_d(y,\lambda)}{G_d(0,\lambda)}.
  \end{align*}
  This yields
  \begin{align*}
      \mathbb{P}(G<\tau_0)=p\sum_{y\neq 0}\left(G_{d}(e_1-y,1-p)-\frac{G_{d}(e_1,1-p)G_{d}(y,\lambda)}{G_{d}(0,1-p)}\right).
  \end{align*}
  We note that for all $x\in\mathbb{Z}^2$,
  \begin{align*}
      \sum_yG_2(x-y,\lambda)=\sum_yG_2(y,\lambda)=G_2(0,\lambda)=\frac{1}{1-\lambda}
  \end{align*}
 and therefore
  \begin{align*}
      \sum_{y\neq 0}G_{2}(x-y,1-p)=\frac{1}{1-(1-p)}-G_{2}(x,1-p)=\frac{1}{p}-G_{2}(x,1-p).
  \end{align*}
  This yields
  \begin{align*}
      \mathbb{P}(G<\tau_0)&=p\left(\frac{1}{p}-G_{2}(e_1,1-p)-\frac{G_{2}(e_1,1-p)}{G_{2}(0,1-p)}\left(\frac{1}{p}-G_{2}(0,1-p)\right)\right)
      \\&=1-\frac{G_{2}(e_1,1-p)}{G_{2}(0,1-p)}.
  \end{align*}
  The proof of the Lemma is finished after combining \ref{a2} and \ref{F2} and applying Lemma \ref{Zrec}. 
  \end{proof}
We continue with the investigation of the return time $Z_1$ in transient dimensions $d\geq 3$. Due to transience, the right quantity to
look at is no longer the Laplace transform of $Z_1$, but instead the probability $\mathbb{P}(Z_1=\infty)$ of never returning to the state $(0,1)$. The following lemma is an analogous version of Lemma \ref{Zrec} for $d\geq 3$:    
\begin{lemma}\label{Ztr}
    For $d\geq 3$ we have
    \begin{align*}
        \mathbb{P}(Z_1=\infty)=1-\frac{2d(\kappa+\rho)^2+s_1\frac{G_d(e_1)G_d(0,q)}{G_d(0,q)-G_d(e_1,q)}}{(s_1+2d(\kappa+\rho))G_d(0)},
    \end{align*}
\end{lemma}
where $G_d(x)$ denotes the Green's function
\begin{align}
    G_d(x):=\mathbb{E}^X_x\left[\int_0^{\infty}\delta_0(X(s))\d s\right]
\end{align}
of a random walk $X$ with jump rate $2d$ in $0$ with start in $x$. 
\begin{proof}
In a similar manner as in the proof of Lemma \ref{Zrec}, using the decomposition \eqref{decomposition}, we observe that
 \begin{align*}
        \mathbb{P}(Z_1=\infty)=\frac{2d(\kappa+\rho)}{s_1+2d(\kappa+\rho)}\mathbb{P}(R_{0,0}^{2d(\kappa+\rho)}=\infty)+\frac{s_1}{s_1+2d(\kappa+\rho)}\mathbb{E}^Y[\mathbb{P}(R_{Y,0}^{2d(\kappa+\rho)}=\infty)],
    \end{align*}
    where $R_{0,0}^{2d(\kappa+\rho)}$ corresponds to the return time of the walk to the origin after starting in origin, as on the event that $Z$ moves away from $0$ (say to $e_1$) before getting dormant, the exponential waiting time of jumping from $0$ to $e_1$ and the return time $R_{e_1,0}^{2d(\kappa+\rho)}$ add up to $R_{0,0}^{2d(\kappa+\rho)}$ (in dimensions $d=1,2$ this exponential jump time was negligible, such that we only considered $R_{e_1,0}^{2d(\kappa+\rho)}$). Now, the probability $\mathbb{P}(R_{0,0}^{2d(\kappa+\rho)}=\infty)$ of never hitting zero is given by
\begin{align*}
   \mathbb{P}(R_{0,0}^{2d(\kappa+\rho)}=\infty)=1-\frac{\kappa+\rho}{G_d(0)}, 
\end{align*}
(cf.\,\cite{lawlerlimic}), as the dependence on the jump rate $2d(\kappa+\rho)$ cancels out in the quotient. Consequently,
\begin{align*}
   \mathbb{P}(Z_1=\infty)&=\frac{2d(\kappa+\rho)}{s_1+2d(\kappa+\rho)}\left(1-\frac{\kappa+\rho}{G_d(0)}\right)+\frac{s_1}{s_1+2d(\kappa+\rho)}\mathbb{E}^Y\left[\left(1-\frac{G_d(Y)}{G_d(0)}\right)\right]
   \\&=1-\frac{2d(\kappa+\rho)^2+s_1\mathbb{E}^Y[G_d(Y)]}{(s_1+2d(\kappa+\rho))G_d(0)}.
\end{align*}
Note that, as $G_d$ is harmonic, Lemma \ref{E[Y]} is applicable again and so
\begin{align*}
    \mathbb{E}[G_d(Y)]=\frac{G_d(e_1)}{\mathbb{P}(G<\tau_0)}=\frac{G_d(e_1)}{1-\frac{G_d(e_1,q)}{G_d(0,q)}}=\frac{G_d(e_1)G_d(0,q)}{G_d(0,q)-G_d(e_1,q)}.
\end{align*}
This finishes the proof of the lemma. 
\end{proof}
\section{Proof of Theorem 1.1}
We are now ready to prove our main result:
\begin{proof}[Proof of Theorem 1.1]
    For $d=1$, by Lemma \ref{d=1} we have that
    \begin{align*}
        1-\lambda\overline{F}(\lambda)=1-\mathbb{E}[e^{-\lambda Z_1}]=\frac{2\sqrt{\lambda}(\kappa+\rho+s_1C_{1,\rho,s_1})}{\sqrt{\kappa+\rho}(s_1+2(\kappa+\rho))}+O(\lambda),\qquad \lambda\to 0.
    \end{align*}
    and 
    \begin{align*}
        1-\mu\lambda\overline{F}(\lambda)=1-\mu+\mu\frac{2\sqrt{\lambda}(\kappa+\rho+s_1C_{1,\rho,s_1})}{\sqrt{\kappa+\rho}(s_1+2(\kappa+\rho))}.
    \end{align*}
    Plugging in $\mu=\frac{2(\kappa+\rho)+s_1}{2(\kappa+\rho)+s_1+\gamma}$ and applying Lemma \ref{G} yields
    \begin{align*}
        \widehat{G}_\mu(\lambda)&=\frac{2(\kappa+\rho+\sqrt{\kappa+\rho}s_1C_{1,\rho,\kappa,s_1})}{\gamma\sqrt{\lambda(\kappa+\rho)}+2\lambda(\kappa+\rho+s_1C_{1,\rho,s_1})}+O(\lambda)
        \\&=\frac{2(\sqrt{\kappa+\rho}+s_1C_{1,\rho,\kappa,s_1})}{\gamma\sqrt{\lambda}}+O(\lambda), \qquad\lambda\to 0.
    \end{align*}
    By a Tauberian theorem we obtain
    \begin{align*}
        \left<U(t)\right>=G_\mu(t)\sim \frac{2(\kappa+\rho+s_1C_{1,\rho,s_1})}{\gamma\sqrt{\pi t(\kappa+\rho)}}, \qquad t\to\infty.
    \end{align*}
    \\In an analogous way, we find that for $d=2$,
    \begin{align*}
        1-\lambda\overline{F}(\lambda)=\frac{4(\kappa+\rho)\pi+s_1C_2}{(s_1+4(\kappa+\rho))\log(1/\lambda)}+O(\lambda),
    \end{align*}
   with $C_2:=\frac{\pi G_{2}(0,q)}{G_2(0,1)+G_{2}(e_1,q)}$, and
    \begin{align*}
        1-\mu\lambda\overline{F}(\lambda)=1-\mu+\mu\frac{4(\kappa+\rho)\pi+s_1C_2}{(s_1+4(\kappa+\rho))\log(1/\lambda)}.
    \end{align*}
    Applying Lemma \ref{G} and $\mu=\frac{4(\kappa+\rho)+s_1}{4(\kappa+\rho)+s_1+\gamma}$, we obtain
    \begin{align*}
        \widehat{G}_{\mu}(\lambda)&\sim\frac{\mu(4\pi(\kappa+\rho)+s_1C_2)}{\lambda((1-\mu)(s_1+4(\kappa+\rho))\log(1/\lambda)+\mu(4(\kappa+\rho)+s_1C_2)}
        \\&\sim\frac{4\pi(\kappa+\rho)+s_1C_2}{\lambda\log(1/\lambda)\gamma},\qquad\lambda\to 0. 
    \end{align*}
    Again, by a Tauberian theorem we obtain
    \begin{align*}
        \left<U(t)\right>\sim \frac{4\pi(\kappa+\rho)+s_1C_2}{\gamma\log(t)}, \qquad t\to 0.
    \end{align*}
Next, we deal with the case $d\geq 3$. Denote by
\begin{align*}
\tilde{G}_d(x,i):=\int_0^{\infty}p_d(x,i,t)\,\d t
\end{align*}
the Green's function of the joint process $(Z,\alpha)$ in $(x,i)$, which has the probabilistic representation
\begin{align*}
\tilde{G}_d(x,i)=\mathbb{E}_{(x,i)}^{(Z,\alpha)}\left[\int_0^{\infty}\delta_{(0,1)}(Z(s),\alpha(s))\,\d s\right].
\end{align*}
Hence, for all $(x,i)\in\mathbb{Z}^d\times\{0,1\}$ the quantity
\begin{align*}
v(x,i):=\lim_{t\to\infty}v(x,i,t)=\mathbb{E}_{(x,i)}^{(Z,\alpha)}\left[\exp\left(-\gamma\int_0^{\infty}\delta_{(0,1)}(Z(s),\alpha(s))\,\d s\right)\right]
\end{align*}
lies in $(0,1)$. Moreover, $v$ solves the boundary value problem
\begin{align*}
\left\{\begin{array}{llll}0&=&(i\kappa+\rho)\Delta v(x,i)-\gamma\delta_{(0,1)}(x,i)v(x,i), &(x,i)\in\mathbb{Z}^d\times\{0,1\},\\[8pt]1&=&\lim_{\|x\|\to\infty}v(x,i), &i\in\{0,1\},\end{array}\right.
\end{align*}
and can therefore be written as
\begin{align*}
v(0,1)=1-\gamma\int_0^{\infty}p_d(0,1,t)v(0,1)\,\d t=1-\gamma v(0,1)\tilde{G}_d(0,1)
\end{align*}
in point $(0,1)$, where $p_d$ denotes the transition density function of $(Z,\alpha)$. Solving for $v(0,1)$ gives
\begin{align}\label{v}
v(0,1)=1-\frac{\gamma \tilde{G}_d(0,1)}{1+\gamma \tilde{G}_d(0,1)}.
\end{align}
The survival probability converges therefore to a non-trivial limit in $(0,1)$ in all dimensions $d\geq 3$. Moreover, note that, $\tilde{G}_d(0,1)$ can be expressed as $1$ over the probability that $Z$ ever comes back to $(0,1)$ again, i.e.
\begin{align}
    \tilde{G}_d(0,1)=\frac{1}{\mathbb{P}(Z_1<\infty)}=\frac{1}{1-\mathbb{P}(Z_1=\infty)}=\frac{(s_1+2d(\kappa+\rho))G_d(0)}{2d(\kappa+\rho)^2+s_1\frac{G_d(e_1)G_d(0,q)}{G_d(0,q)-G_d(e_1,q)}},
\end{align}
by using Lemma \ref{Ztr} (observe that we obtain $\tilde{G}_d(0,1)=G_d(0)$ for $s_1=0$, such that this is consistent with the case without dormancy). Plugging this into \eqref{v} yields
\begin{align*}
    v(0,1)&=1-\frac{\gamma \left(G_d(0)+\frac{s_1}{2d(\kappa+\rho)}\right)}{\kappa+\rho +\gamma \left(G_d(0)+\frac{s_1K_d}{2d(\kappa+\rho)}\right)}
\end{align*}
with $K_d:=K_{d,s_1,\rho}$ defined as
\begin{align*}
K_d:=\frac{G_d(e_1)G_d(0,q)}{G_d(0,q)-G_d(e_1,q)}.
\end{align*}
\end{proof}
\appendix
\section{Some auxiliary results}
In this appendix we prove two results which we used in the proof of Lemma \ref{d=2}:
\begin{lemma}\label{lem:error_E_goes0}
Abbreviate $\phi(k):=\frac{1}{2}(\cos(k_1)+\cos(k_2))$ for $k=(k_1,k_2)\in[-\pi,\pi]^2$ and let
\begin{align*}
    E(x,\lambda)&:=\int_{[-\pi,\pi]^2}\left(\frac{1-e^{ik\cdot x}}{1+\lambda-\phi(k)}-\frac{1-e^{ik\cdot x}}{1-\phi(k)}\right)\frac{\emph{d} k}{(2\pi)^2}.
\end{align*}
Then $\mathbb{E}[E(W,\lambda)]\rightarrow 0$ as $\lambda\rightarrow 0$ for any random variable $W$ with finite variance.
\end{lemma}
\begin{proof}[Proof]
    Note that for every $\delta>0$,
    \begin{align*}
        E(x,\lambda)&=\int_{[-\pi,\pi]^2}\frac{-\lambda(1-e^{ik\cdot x})}{(1+\lambda-\phi(k))(1-\phi(k))}\frac{\d k}{(2\pi)^2}\\
        &=\int_{[-\pi,\pi]^2\setminus B_\delta(0)}\frac{-\lambda(1-e^{ik\cdot x})}{(1+\lambda-\phi(k))(1-\phi(k))}\frac{\d k}{(2\pi)^2}
        \\&+\int_{B_\delta(0)}\frac{-\lambda(1-e^{ik\cdot x}+ik\cdot x)}{(1+\lambda-\phi(k))(1-\phi(k))}\frac{\d k}{(2\pi)^2}\\
        &=:K_\lambda(\delta,x)+I_\lambda(\delta,x),
    \end{align*}
where in the second line we considered what happens in a ball of radius $\delta$ around $0$, as well as outside of that separately, as the point $0$ is where there is a singularity in the integral. We also added and subtracted $ik\cdot x$ in the second term $I_\lambda(\delta,x)$. Observe that for all $\delta>0$,
    \begin{align*}
        K_\lambda(\delta,x)\leq C_1\left(\frac{1}{1-\cos(\delta)}\right)^2\lambda,
    \end{align*}
    for some constant $C_1>0$ independent of $\lambda$ and $\delta$. To deal with $I_\lambda(\delta,x)$ we note that
    \begin{align*}
        \frac{\lambda}{\left|1+\lambda-\phi(k)\right|}\leq1
    \end{align*}
    and $(1-e^{ik\cdot x}+ik\cdot x)\leq |k|^2|x|^2$ as well as $\cos(x)\leq1-\frac{1}{4}x^2$ for sufficiently small $x$, which implies $1-\phi(k)\geq\frac{1}{8}(k_1^2+k_2^2)$.
    Therefore, 
    \begin{align*}
        \left|\frac{1-e^{ik\cdot x}+ik\cdot x}{1-\phi(k)}\right|\leq C_2|x|^2
    \end{align*}
    for all $(k_1,k_2)\in B_\delta(0)$ and some constant $C_2>0$, and thus, together with the last inequality, we have that $I_\lambda(\delta,x)\leq C_2|(B_\delta(0)||x|^2\leq C_3\delta^2|x|^2$ for some other constant $C_3>0$. This yields
    \begin{align*}
        E(x,\lambda)=K_\lambda(\delta,x)+I_\lambda(\delta,x)\leq C_4\left(\lambda\left(\frac{1}{1-\cos(\delta)}\right)^2+\delta^2|x|^2\right)
    \end{align*}
    for some absolute constant $C_4>0$, which is valid for all sufficiently small $\delta$. Taking the expectation with respect to $W$ and noting that $\mathbb{E}[W^2]<\infty$ we get,
    \begin{align*}
        \limsup_{\lambda\rightarrow0}\mathbb{E}[E(W,\lambda)]\leq C_4\delta^2\mathbb{E}[W^2].
    \end{align*}
    Since this holds for all sufficiently small $\delta$, taking $\delta\rightarrow0$ gives the result.
\end{proof}
Next, we define the \emph{Green's generating function} of some point $x\in\mathbb{Z}$ as
\begin{align*}
    G_2(x,\lambda):=\int_{[-\pi,\pi]^2}\frac{e^{ik\cdot x}}{1+\lambda-\phi(k)}\frac{dk}{(2\pi)^2},
\end{align*}
where the subscript $2$ refers to the dimension. The next lemma provides asymptotics for the Green's generating function at $x=0$:
\begin{lemma}\label{lem:asymptotics_G(0)}
    We have that
    \begin{align*}
        G_2(0,\lambda)=\frac{1}{\pi}\log(1/\lambda)+O(1), \qquad \lambda\to 0.
    \end{align*}
\end{lemma}
\begin{proof}[Proof]
    For $\delta>0$ small enough we can decompose
    \begin{align*}
        G_2(0,\lambda)&=\int_{[-\pi,\pi]^2\setminus B_\delta(0)}\frac{1}{\lambda+1-\phi(k)}\frac{\d k}{(2\pi)^2}+\int_{B_\delta(0)}\frac{1}{\lambda+\frac{1}{4}(k_1^2+k_2^2)}\frac{\d k}{(2\pi)^2}\\
        &+\int_{B_\delta(0)}\frac{1}{\lambda+1-\phi(k)}-\frac{1}{\lambda+\frac{1}{4}(k_1^2+k_2^2)}\frac{\d k}{(2\pi)^2}
    \end{align*}
    We will show that for a fixed $\delta>0$, only the middle term will contribute to the asymptotics, as $\lambda\rightarrow0$, and that the other terms will be of constant order. Concerning the first term, we observe that the integrand and domain of integration are bounded as $\lambda\rightarrow0$. The last term can be rewritten as
    \begin{align}\label{last}
        -\frac{1}{2}\int_{B_\delta(0)}\frac{(1-\cos(k_1)-\frac{1}{2}k_1^2)+(1-\cos(k_2)-\frac{1}{2}k_2^2)}{(\lambda+1-\phi(k))(\lambda+\frac{1}{4}(k_1^2+k_2^2))}\frac{\d k}{(2\pi)^2},
    \end{align}
and by a Taylor approximation we have that $|1-\cos(k_1)-\frac{1}{2}k_1^2|\leq\frac{1}{24}k_1^4$. Since we chose $\delta>0$ small enough we have $1-\cos(x)\geq\frac{x^2}{4}$ for $x\in B_\delta(0)$ and so \eqref{last} is bounded by
    \begin{align*}
        C\int_{B_\delta(0)}\frac{k_1^4+k_2^4}{(\lambda+\frac{1}{8}(k_1^2+k_2^2))(\lambda+\frac{1}{4}(k_1^2+k_2^2))}\frac{\d k}{(2\pi)^2}
    \end{align*}
    for some absolute constant $C>0$. This last expression is bounded as $\lambda\rightarrow0$. Next, a computation in polar coordinates yields
    \begin{align*}
        \int_{B_\delta(0)}\frac{1}{\lambda+\frac{1}{4}(k_1^2+k_2^2)}\frac{\d k}{(2\pi)^2}&=\frac{1}{(2\pi)^2}\int_0^{2\pi}\int_0^\delta\frac{r}{\lambda+\frac{r^2}{4}}\d r\d\theta\\
        &=\frac{1}{\pi}\log(1/\lambda)+\frac{1}{\pi}\log(\lambda+\frac{\delta^2}{4}).
    \end{align*}
    Since $\delta>0$ is fixed, the term $\frac{1}{\pi}\log(\lambda+\frac{\delta^2}{4})$ does not contribute to the asymptotics, which completes the proof of our claim. 
\end{proof}

\subsection*{Acknowledgement}
The authors would like to thank Professor Alison Etheridge for her insightful comments.


\begin{thebibliography}{WK2020}


\bibitem[BYZ13]{baran}
\newblock {\sc N.~A.~Baran, G.~Yin, C.~Zhu},
\newblock {\em Feynman-Kac formula for switching
diffusions: connections of systems of partial
differential equations and stochastic
differential equations},
\newblock Advances in Difference Equations, No. 315 (2013).

\smallskip

\bibitem[BHS21]{blath}
\newblock {\sc J.~Blath, F.~Hermann, M.~Slowik},
\newblock {\em A branching process model for dormancy and seed banks in
randomly fluctuating environments},
\newblock Journal of Mathematical Biology 83, No. 17 (2021).

\smallskip


\bibitem[BHLWB21]{dormancy}
\newblock {\sc J.~Blath, J.~T.~Lennon, F.~den Hollander, M.~Wilke Berenguer},
\newblock {\em Principles of seed banks and the emergence of
complexity from dormancy},
\newblock Nature Communications 12, No. 4807 (2021).

\smallskip

\bibitem[C66]{cohen}
\newblock {\sc D.~Cohen},
\newblock {\em Optimizing reproduction in a randomly varying environment},
\newblock Journal of Theoretical Biology 16, 267–282 (1966).

\smallskip

\bibitem[GM90]{garmol}
\newblock {\sc J.~Gärtner, S.~Molchanov},
\newblock {\em Parabolic problems for the Anderson
model I. Intermittency and related topics},
\newblock Communications in Mathematical Physics 132, 613–
655 (1990).

\smallskip 

\bibitem[HJV07]{bookhaccou}
\newblock {\sc P.~Jagers, P.~Haccou, V.~A.~Vatutin},
\newblock {\em Branching processes: variation, growth, and extinction of populations},
\newblock Biometrics 62, No. 4, 1269–1270 (December 2006).

\smallskip 

\bibitem[K16]{PAM}
\newblock {\sc  W.~König}, 
\newblock {\em The Parabolic Anderson Model: Random Walk in Random Potential}, 
\newblock Birkhäuser (2016).

\smallskip

\bibitem[L96]{lawler}
\newblock {\sc G.~F.~Lawler},
\newblock {\em Intersections of random walks},
\newblock Birkhäuser Boston (1996).


\smallskip

\bibitem[LL10]{lawlerlimic}
\newblock {\sc  G.~F.~Lwaler, V.~Limic}, 
\newblock {\em Random walk: A modern introduction}, 
\newblock Cambridge University Press (2010).

\smallskip

\bibitem[SW11]{sw}
\newblock {\sc A.~Schnitzler, T.~Wolff},
\newblock {\em Precise asymptotics for the parabolic Anderson
model with a moving catalyst or trap},
\newblock Conference Paper. In {\em Probability in complex physical systems}, 69-89, Springer (2012)

\smallskip

\bibitem[S24]{shafigh}
\newblock {\sc H.~Shafigh},
\newblock {\em A spatial model for dormancy in random environment},
\newblock ArXiv e-prints (September 2024), arXiv:2304.07703  

\smallskip

\bibitem[S25]{shafigh25}
\newblock {\sc H.~Shafigh},
\newblock {\em Dormany in random environment: Symmetric exclusion},
\newblock ArXiv e-prints (Januar 2025), arXiv:2501.03807

\smallskip

\bibitem[YZ10]{switching}
\newblock {\sc G.~Yin, C.~Zhu},
\newblock {\em Hybrid Switching Diffusions: Properties and Applications}
\newblock Springer New York, NY (2010)


\end{thebibliography}
\end{document}